\documentclass{jancl}
\usepackage{graphicx}
\usepackage{amssymb,amsmath,makeidx,euscript}

\title[Generalized Kracht formulas]
{An extension of Kracht's theorem\\ to generalized Sahlqvist formulas}
\journal{0}{0}{0}{-}{}
\author{Stanislav Kikot}
\address{Moscow State University\\
Moscow, Vorobjevy Gory, 1, (Russia),\\[3pt] staskikotx@gmail.com}
\abstract{
Sahlqvist formulas are  a syntactically specified class
of modal formulas proposed by Hendrik Sahlqvist in 1975. They are 
important because of their first-order definability and 
canonicity, and hence axiomatize complete modal logics. The
first-order properties definable by Sahlqvist formulas were
syntactically characterized by Marcus Kracht in 1993. 
The present paper extends Kracht's theorem to the class of `generalized
Sahlqvist formulas' introduced by Goranko and Vakarelov and describes an
appropriate generalization of Kracht formulas.
}

\def\iffa{\Longleftrightarrow}
\def\d{\Diamond}
\def\Ml{{\cal ML}}
\def\b{\Box}
\def\K{{\EuScript K}}
\def\E{{\cal E}}
\def\P{{\cal P}}

\def\P{{\cal P}}
\def\rank{\mbox{ rank }}
\def\tr{\triangleright}

\begin{document}
\maketitle
\section{Introduction}

The Sahlqvist theorem is a hard working horse  in modal logic. 
It describes a large class of first-order definable canonical modal
formulas.
A standard proof of completeness results boils down to
finding relevant first-order properties and corresponding Sahlqvist formulas and
next --- to  applying
Sahlqvist completeness theorem.  Also Sahlqvist formulas are often applied
for proofs of negative results such as non-finite axiomatizability.

Kracht's theorem is an important addition to the Sahlqvist theorem. It 
explicitly describes the class of first-order correspondents to Sahlqvist
formulas \cite{Kr1}, \cite{Kr2}. Moreover, it gives an algorithm constructing
a Sahlqvist formula from its first-order analogue.

So when we encode a first-order condition 
into a Sahlqvist formula, we implicitly use Kracht's algorithm.  That is
why for axiomatizing modal logics 
Kracht's theorem is not less important than
the Sahlqvist theorem.

In \cite{Gor1}, \cite{Gor2} the Sahlqvist theorem was further generalized.
These results turned out to be at the intersection 
of at least two  research lines.

The first line came from attempts at axiomatizing 
many-dimensional modal logics.
Probably,  the first known generalized Sahlqvist formula was 
$cub_1$ (see page \pageref{cub1} of this paper) 
for the first time published in  \cite{Sh78}, 
expressing the `cubifying' property of 3-dimensional
product frames (see Figure 1)
\begin{equation}\label{cub_fo}
\begin{array}{l}
(\forall x_0) \forall x_1 \forall x_2 \forall x_3 \left( x_0 R_1 x_1 \land
x_0 R_2 x_2 \land
x_0 R_3 x_3 \to
\right.\\
\to
\exists y
(
y \in  R_3(R_2(x_1)\cap R_1(x_2))\land y\in  R_2(R_3(x_1)\cap R_1(x_3))\land
\\
\left.\land y \in R_1(R_2(x_3)\cap R_3(x_2)))\right).\\
\end{array}
\end{equation}
  Modifications of this formula were used by 
A. Kurucz in the  proof of  some
negative results on $\ge 3$-dimensional products \cite{Agi1}, 
\cite{Agi2}. Let us also mention that
generalized Sahlqvist formulas appear in axiomatizing
2-dimensional squares with extinguished diagonal \cite{K3}.

\begin{figure}[cub]{}
\centerline{\includegraphics{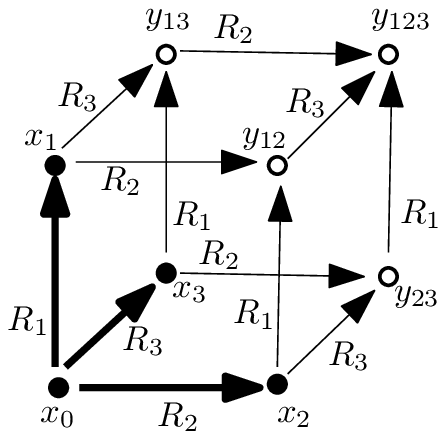}}
\end{figure}

First-order conditions like (\ref{cub_fo}) can be illustrated by
pictures with black and white points, and bold and simple arrows, as
in Figure 1. A formal analogue of such a picture is  the notion of  
a {\sl diagram}.  It turns out
that under some natural conditions, the corresponding 
first-order $\forall\exists$-formula is modally definable 
if and only if the diagram does not have non-oriented cycles consisting of 
white points and  simple arrows, and, in the case of modal
definability, the $\forall\exists$-formula always 
corresponds to a generalized Sahlqvist formula   \cite{K1}.

The second line of research arises from the natural problem --- to find
sufficient conditions for 
first-order definability and canonicity of modal formulas. 
The relevant part starts with
\cite{Gor1} extending the Sahlqvist theorem to polyadic modal
languages. 
The same paper gives an example of  first-order definable 
and canonical modal formulas that are not Sahlqvist (namely, 
formula $D_2$ from Example \ref{example} below). However,
the question if these new formulas have 
Sahlqvist equivalents, was remaining unsolved for some time. 
This question was solved by V. Goranko and D. Vakarelov who introduced
the notion of a-persistence and showed that 
all Sahlqvist formulas are a-persistent while $D_2$ is not \cite{Gor2}. 
It was in \cite{Gor2} that the notion of a `generalized Sahlqvist formula', 
that lies in the center of the present paper,  was introduced 
as a partial case of so called `inductive formula'.

Then algorithms were proposed 
(see \cite{CGV1}, \cite{CGV2} and references therein), 
for computing first-order equivalents of some
modal formulas. The Sahlqvist theorem was further generalized 
in \cite{V2} and \cite{V}, 
yielding the class
of  complex Sahlqvist formulas,
but  they are actually semantically equivalent to 
standard Sahlqvist formulas.

Another challenging problem is: 
`given a first-order formula, find the modal logic of the corresponding
elementary class'. Let us mention that the problems
`given a first-order formula, determine if it is modally definable' and
`given a modal formula, determine if it is first-order definable' are
undecidable due to Chagrova's theorem \cite{ChL}. That is why any sufficient condition
for modal (or f. o.) definability is very interesting by itself.
In this context besides the above cited Kracht's result \cite{Kr1}, \cite{Kr2} 
and the study of diagram formulas \cite{K1} we can make especially mention   
the brilliant work \cite{Hodkinson} giving an explicit infinite axiomatization
for any elementary class. In some particular cases more
concise (although also infinite) axiomatizations are constructed in \cite{BSS}.
However, we still do not have a criterion of finite axiomatizability for 
the logics from  \cite{Hodkinson} and \cite{BSS}.

The present paper continues the study on modal logics of 
elementary classes. We extend the
class of Kracht formulas to the class of `generalized Kracht formulas'.
Then we propose an algorithm  constructing a modal correspondent for
a given generalized Kracht formula. This modal correspondent is a generalized
Sahlqvist formula, and therefore it is canonical (and, a fortiori, Kripke complete).

Our terminology slightly differs from \cite{Gor2}; in particular,
the notion `regular formula' has a different  meaning.   
Also, the term `safe expression' is not the same as in  \cite{Bl}.

\section{Regular box-formulas.}
We consider the modal language $\Ml_\Lambda$ with countably many propositional variables,
unary modalities $\d_\lambda$ and
their duals $\b_\lambda$, where $\lambda\in \Lambda$,  boolean connectives $\land,
\lor, \neg, \to $ and boolean constants $\top, \bot$. A formula in this language
is called {\sl positive} if it does not contain $\neg$ and $\to$ (but
may contain $\bot$.)

Recall that in the Sahlqvist theorem  `boxed atoms' 
(i.e. the expressions of the form $\b^n p$) are
crucial,
because they allow us to obtain the minimal valuation for an antecedent. 
In generalized Sahlqvist formulas  `boxed atoms' are replaced by
`regular box-formulas'.

\begin{definitionCite}{Gor2} \label{defregvak}
{\sl A box-formula } is  defined by recursion:
\begin{itemize}
\item a variable  $p_i$ is a box-formula; 
\item if $POS$ is a positive modal formula and $BF$ is a box-formula 
then $POS\to BF$ is a box-formula;
\item if $BF$ is a box-formula then $\b_\lambda BF$ is a box-formula.
\end{itemize}  
 
Thus a box-formula is equivalent  to one of the form 
$$POS_1 \to \b^{\alpha_1}(POS_2 \to  \b^{\alpha_2}(POS_3 \to \ldots \to p_i)\ldots),$$
where $\b^{\alpha_j}$ are sequences of boxes, $POS_j$ are positive.

The last variable $p_i$  of this formula is called its {\sl head}.
$BF\succ p_i$ denotes that $p_i$ is the head of a box-formula $BF$.

Let $A$ be a set of box-formulas. The {\sl dependency graph } of 
$A$ is an oriented graph $G=(V_A, E_A)$, 
where the set of vertices $V_A$ = $\{p_1, \ldots, p_n\}$ consists of
all variables occurring in $A$, and the adjacency relation is 
$$p_i E_A p_j \iffa p_i \mbox{ occurs (not as a head) in some formula } 
\phi\in A \mbox{ with the head } p_j.$$ 
A set of box-formulas $A$ is called {\sl regular} if its dependency graph is acyclic, \ie, it does not contain oriented cycles.
\end{definitionCite}

%Strictly speaking, in \cite{Gor2} the notion of  a regular set of box-formulas is not %defined explicitly. Instead of this 
%the dependency graph of all box-formulas used in the 
%antecedent of a generalized Sahlvist formula was considered

We will use a more convenient technical version of Definition \ref{defregvak}.

The set of propositional variables is split into
countably many groups $p^0_1, p^0_2, p^0_3, \ldots$,
$p^1_1, p^1_2, p^1_3, \ldots$, $p^2_1, p^2_2, p^2_3, \ldots$ and so on.
The upper index (called the  {\sl rank})  is the number of the group and the lower index is
the number of a variable within a group. Put $\bar p^i=
\{p^i_1, p^i_2, p^i_3, \ldots\}$.

\begin{definition}\label{defregmy}
{\sl A regular box-formula of rank $k$} is defined by recursion:
\begin{itemize}
\item a variable $p^k_i$ is a box-formula of rank $k$,
\item if $POS(\bar p^0,\bar p^1, \ldots, \bar p^{k-1})$ is a positive modal formula,
depending only on the variables of rank $<k$ and $REG$ is a 
regular box-formula of rank $k$ then $POS(\bar p^0, \bar p^1, \ldots, \bar p^{k-1})\to REG$ is a regular box-formula of rank $k$,
\item if $REG$ is a 
regular box-formula of rank $k$ then $\b_\lambda REG$ is a regular 
box-formula of rank $k$.
\end{itemize}  
\end{definition}

\begin{lemma}
Let $A$ be a set of modal formulas. Then 

(1) if $A$ is a set of regular box-formulas (in the sense of Definition \ref{defregmy}), then
$A$ is a regular set of box-formulas (in the sense of Definition \ref{defregvak}).

(2) if $A$ is a regular set of box-formulas, then
we can range the propositional variables (\ie  choose the  upper indices) 
so  that 
$A$ becomes a set of regular box-formulas.
\end{lemma}

\begin{proof*} 
(1) is trivial. In fact, if $A$ is a set of regular formulas and $p^s_i E_A p^t_j$, then $s<t$.  So the dependency graph 
does not contain oriented cycles.

We prove (2) by induction on the number of vertices in $V_A$.
If it has a single vertex, the statement is trivial. Suppose it has $n$ vertices.
Since  our graph does not have oriented cycles, there is a vertex $v$ in $V_A$ without successors. Suppose $v$ corresponds to a variable $p_l$ for some
$l \le n$. We eliminate this vertex (and of course, all entering edges) 
and obtain the graph $G'_A=(V_A', E_A')$.
Since $v$ does not have
successors, $p_l$ can occur only in the heads of box-formulas from $A$.
If $A'$ is obtained from $A$ by eliminating box-formulas with the head $p_l$, then
$G_{A'}=G'_A$.
By the induction hypothesis we can range the vertices of  $G'_A$ so that
all formulas in $A'$ become regular. For $i \neq l$ let  $r(i)$ be 
the rank of $p_i$. Put the rank of 
$p_l$ to be $\max r(i)+1$. Then $A$ is a set of regular formulas, since all formulas in
$A\setminus A'$ have $p_l$ as their head, and the rank of $p_l$ is maximal.
\end {proof*}

%We introduce some more notation to deal with the minimal valuation. 
%In fact, our computation is similar to \cite{Gor2}, 
%but our notation yields  a more explicit form of
%a minimal valuation. The
%operators $KV$ and $KF$ described below 
%correspond to the minimal valuation for a fixed variable for the cases 
%where the valuations for the variables of smaller rank are
%respectively fixed and arbitrary. The operators $KP$ and $KVF$ are
%auxiliary, we need them only  for formal description of $KV$ and $KF$. 

%The rest of this Section is rather technical computation of exact formulas
%for minimal valuations.  The reader not interested 
%in the proofs might proceed to Sections 3 and 5.

Besides  the modal language $\Ml_\Lambda$, we  need additional  languages 
$L^\#_k$,$L^P_k$ and $L$. Their vocabularies are
\begin{itemize}
\item for $L^P_k:   P^l_i \,(l<k), \cap, \cup, R_\lambda^{-1}, R_\lambda^{\b},\top,\bot;$
\item for $L^\#_k:   \#, P^l_i\, (l<k), \cap, \cup, R^{-1}_\lambda, R_\lambda^{\b}, R_\lambda,\top,\bot;$
\item for $L:   x_i, \cap, \cup, R^{-1}_\lambda, R^{\b}_\lambda, R_\lambda,\top,\bot.$
\end{itemize}

Here $\bot, \top$ are constants, $P^l_i, \#, x_i$ are variables,
$R_\lambda^{-1}, R_\lambda^{\b}, R_\lambda$ are unary function symbols,
$\cap, \cup$ are binary function symbols. 
We call the terms of these languages  {\sl  expressions}.

To every  regular box-formula $\phi$ of rank $k$ we assign an  $L^\#_k$-expression $KV^\phi$. (Later we shall see that $KV^{\phi}$ is the operator for the relative minimal valuation for the head of $\phi$.)

First we assign an expression $KP^{POS} \in L^P_k$ to every positive formula $POS$ :
$$KP^{\top} := \top,\quad KP^{\bot} := \bot,$$ 
$$KP^{p^l_i} := P^l_i, \mbox{ where } l<k,$$ 
$$KP^{POS_1\land POS_2} := KP^{POS_1}\cap KP^{POS_2},$$ 
$$KP^{POS_1\lor POS_2} := KP^{POS_1}\cup KP^{POS_2},$$  
$$KP^{\d_\lambda POS} := R^{-1}_\lambda(KP^{POS}),$$ 
$$KP^{\b_\lambda POS} := R^{\b}_\lambda(KP^{POS}).$$

This definition obviously corresponds to the truth definition in the
standard Kripke semantics. If we have a frame
$F=(W,(R_\lambda:\lambda\in \Lambda))$ and  $\theta$ is a
valuation for the variables $p^l_i$, where $l<k$,
then $\theta(POS)$ is the value of $KP^{POS}$
under the interpretation  $I$ sending $\top$ to $W$, $\bot$ to $\emptyset$,  
$P^l_i$ to $\theta(p^l_i)$,
$R^{-1}_\lambda(A)$ to 
$\{ x\mid \exists y \, x R y \mbox{ and } y \in A\}$, $R^{\b}_\lambda(A)$ to
$\{ x\mid \forall y \mbox{ if } x R y \mbox{ then } y \in A\}$

Now we assign an  $L^\#_k$-expression $KV^\phi$ to any regular box-formula $\phi$ of rank $k$. 

\begin{definition}\label{K^phi}
We set $$KV^{p^k_i}: = \#,$$
$$KV^{POS \to \psi}: = KV^\psi(\#\cap KP^{POS}),$$
 $$KV^{\b_\lambda\psi}: =  KV^\psi(R_\lambda(\#)).$$
\end{definition}

Here $KV^\phi(t)$
denotes the substitution instance $[t/\#] KV^\phi$.
That is to obtain $KV^{POS \to \psi}$, we substitute 
the term $\#\cap KP^{POS}$ for $\#$ in $KV^\psi$, 
and to obtain
$KV^{\b_\lambda \psi}$, we substitute
the term $R_\lambda(\#)$ for  $\#$ in  $KV^\psi$.

\begin{example}
\begin{enumerate}
\item Let $\phi = \b^l_\lambda p^1_0$. Then $KV^\phi= R^l_\lambda(\#)$.
\item If $\phi = \b_1 (\d_2 p_0^0\to\b_3 p^1_0)$, then $KV^\phi 
  = R_3(R_2^{-1}(P^0_0)\cap R_1(\#))$.
\end{enumerate}
\end{example}

%In  a frame $F=(W,R)$, we can interpret $KV^\phi$ as
%an operator $KV^\phi : \P(W)^n \times \P(W) \to \P(W)$ in
%the boolean algebra of $F$. 
%ere $n$ denotes the number of different $P^l_i$'s in $\phi$. We mark out 
%the last argument of  $KV^\phi$ and substitute it for $\#$. 

In a model $M=(W, R_\lambda, \theta)$, where $x\in W$,
we can evaluate $KV^\phi(x)$ under the interpretation $I$ described above
and identify it with a certain subset of $W$. 

\begin{lemma}[on monotonicity of $KV^\phi$]\label{mon}
$KV^\phi(x)$ is monotonic with respect to $P^l_i$.
\end{lemma}
\begin{proof*} This is trivial, since all operations
$\cap, \cup, R^{-1}_\lambda, R^{\b}_\lambda, R_\lambda$ are monotonic.
\end{proof*}

The next lemma shows that the operator $KV^\phi$ 
really defines the  `relative minimal valuation' for the truth of $\phi$
in the standard Kripke semantics.

\begin{lemma}[on minimality of $KV^\phi$]\label{min} 
Let $\phi$ be a regular box-formula with a head $p^k_i$.
Consider a Kripke model $M=(W, (R_\lambda:\lambda\in\Lambda), 
\theta)$ where $\theta(p^l_i)=P^l_i$ $(l\le k)$.
Then $$M,x \models\phi \iffa P^k_i \supseteq KV^\phi(x).$$
\end{lemma}

\begin{proof*} The proof is by induction on the length of $\phi$. 
If $\phi$ is a variable, there is nothing to prove. 

Let $\phi = POS \to \psi$. Then

$$x\models \phi \iffa x\models POS \to \psi \iffa$$ 
$$(\mbox{ if }x\models POS, \mbox{ then } x\models\psi)\iffa$$
$$(\mbox{ if }x\models POS, \mbox{ then } P^k_i \supseteq KV^\psi(x)) \iffa$$
$$(x\not\models POS \mbox{ or } P^k_i \supseteq KV^\psi(x))\iffa$$
\begin{equation}\label{disj}
\{x\}\cap KP^{POS} =\emptyset \mbox{ or } P^k_i \supseteq KV^\psi(x)
\end{equation}

There are only two possible values of $\{x\} \cap KP^{POS}$, viz.  
$\{x\}$ and $\emptyset$.
A simple induction argument shows that $KV^{\psi}(\emptyset)=\emptyset$.  
So by an easy study of cases (\ref{disj}) is equivalent to 
$$P^k_i \supseteq KV^\psi(\{x\}\cap KP^{POS}) \iffa$$
$$P^k_i \supseteq KV^\phi(x).$$ 

Let $\phi = \b_\lambda\psi$. 
$$x\models \phi \iffa x\models \b_\lambda\psi \iffa$$ 
$$\forall y (xR_\lambda y \Rightarrow y\models \psi) \iffa$$
$$\forall y (xR_\lambda y \Rightarrow P^k_i\supseteq KV^\psi(y)) \iffa
\footnote{Here we use the fact that $KV^\psi$
is destributive over arbitrary unions.}$$
$$P^k_i \supseteq  KV^\psi(R_\lambda(x)) \iffa$$
$$P^k_i \supseteq KV^\phi(x).$$ 
\end{proof*}

Let $A$ be a finite set of  regular box-formulas, $\P(A)$ be the set of all
subsets of $A$.

Consider a set $V=\{x_1, \ldots, x_n\}$, and a function $f: V \to \P(A)$.

\begin{definition} \label{K_f}
To every variable $p^k_i$  we assign
an expression $KF^{p^k_i}_f$ of our language $L$ (see above) by induction on $k$.
 We put
 $$KF_f^{p^k_i} =  \bigcup_{\phi \succ p^k_i, \phi \in f(x_i)}
KVF^\phi_f(x_i),$$
where $KVF^\phi_f(x_i)$ is obtained by substituting expressions
$KF_f^{p_i^l}$  for $P^l_i$ for all $l<k$ in the expression $KV^\phi(x_i)$,
that we can denote by
$$KVF^\phi_f(x_i) = [KV_f^{p^l_i} / P^l_i]_{l<k} KV^\phi(x_i).$$
In particular, 
$$KF^{p^0_i}_f = \bigcup_{\phi \succ p^0_i, \phi \in f(x_i)}
KV^\phi_f(x_i),$$ where $KV^\phi_f(x_i)$ does not
contain $P$'s.  
\end{definition}

\begin{example} If $V=\{x_1, x_2\}$, $f(x_1)=\{\b_4 p^0_0\}$ and
$f(x_2) = \b_1 (\d_2 p_0^0\to\b_3 p^1_0)$, then
$$KF_f^{p^0_0} = R_4(x_1),$$
$$KF_f^{p^1_0} = R_3\left(R_2^{-1}(R_4(x_1))\cap R_1(x_2)\right).$$
\end{example}

The next lemma shows that the operator $KF^{p^k_i}_f$ corresponds to  the absolute minimal valuation
for a variable $p^k_i$.

\begin{lemma}\label{theta_min}
%In a frame $F=(W, R_\lambda)$ with points  $x_j\in W$ 
Among all valuations $\theta$ such that for all $j$ $x_j \models f(x_j)$
\footnote{Strictly speaking, in this lemma we mean that we have a  frame
$F=(W, (R_\lambda:\lambda\in\Lambda))$ and a valuation of
object variables $g: V \to  W$, so this formula must be read as $F, g(x_j), \theta \models f(x_j)$,
but following Kracht \cite{Kr2} we will identify $x_i$ with $g(x_i)$,
and will not take care of the frame $F$.}
there is the smallest one $\theta_{\min}$, and $\theta_{\min}(p_i^l)=KF_f^{p_i^l}$.
\end{lemma}

\begin{proof*} 
Put $$\rank (f) = \max_{
\begin{array}c
x_j\in V\\ \phi\in f(x_j)
\end{array}} \rank (\phi),$$ where $\rank(\phi)$ denotes the rank of its head. 
Let us introduce a new function $f^-: V\to A$,  as
follows:
$$f^-(x_j)=f(x_j)\cap\{\phi\mid \rank(\phi)<\rank(f)\}.$$
It is clear that 
$$\rank(f^-)\le\rank(f)-1$$

We argue by induction on $\rank f$.

The base: $\rank f =0$.
Then 
$$\theta_{\min}(p^0_i)=\bigcup_{\phi\succ p^0_i, \phi \in f(x_j)} KV^\phi(x_j)=KF_f^{p^0_i}.$$

The induction step. Suppose $\rank f =k.$ 
Consider the map  $f^{-}$. Then by the induction hypothesis there exists
$\theta^-_{\min}(p^l_i)=KF_{f^-}^{p^l_i}$ for $l<k$, such that for any valuation $\theta^-$, given on the variables of $\rank < k$
$$\forall j\  \theta^-, x_j \models f^-(x_j) \to 
\theta^- \supseteq \theta^-_{\min}.$$ Put 
 $$\theta_{\min}(p^k_i) =  KF^{p^k_i}_f=\bigcup_{\phi \succ p^k_i, \phi \in f(x_j)}
KVF_f^\phi(x_j).$$
Suppose that for some $\theta$
$$\forall j\  F,  x_j, \theta \models f(x_j).$$
Let us prove that
$$\theta\supseteq\theta_{\min}.$$
Let $\theta^-$ be a restriction of $\theta$ to  variables of rank  $ < k$.
By the induction hypothesis 
$$\theta^{-}\supseteq\theta_{\min}^{-}.$$
Consider  an arbitrary   $\phi \in f(x_j)$
with the head $p^k_i$. 
By Lemma \ref{min} (on the minimality of $KV^\phi$)
$$\theta(p^k_i) \supseteq KV^\phi(x_j)$$ and by Lemma \ref{mon} (on the monotonicity of $KV^\phi$)
$$KV^\phi(x_j)\supseteq KVF_f^\phi(x_j),$$ hence
$$\theta(p^k_i) \supseteq KVF_f^\phi(x_j).$$ So
$$\theta(p^k_i)\supseteq\bigcup_{\phi \succ p^k_i, \phi \in f(x_j)}  KVF_f^\phi(x_j)=\theta_{\min}(p^k_i)$$
\end{proof*}

\section{Safe expressions}
In this section we study the  values of $KVF$ and $KF$.
\begin{definition} Let $B$ be a set of $L$-expressions.
A {\sl positive combination of $B$} (denoted POS(B)) is any $L$-expression
built from the members of $B$ using only $\cap, \cup, R^{-1}_\lambda, R^{\b}_\lambda, \top, \bot$ (i.~e. all operations of $L$ excepting $R_\lambda$).
\end{definition}

\begin{definition}
Let $\K$ be  the minimal class of $L$-expressions satisfying the conditions:
\begin{itemize}
\item $\{ x_1, \ldots, x_n\} \subseteq \K$,
\item if $S \in \K$, then $R_\lambda(S) \in \K$,
\item if $B \subseteq \K$ and $S\in \K$ then $S\cap POS(B)\in\K$,
\end{itemize}
\end{definition}
where  POS(B) denotes any positive combination of $B$.

Now we give another description of $\K$.

\begin{definition}\label{safeexpr}
Let $\psi$ be a subexpression of $\phi\in L$.
We say that a subexpression $\psi$ is {\sl safe for $\phi$} 
if one of the following holds:

1) $\psi=x_i$;

2) $\psi = R_\lambda(\psi ')$, where $\psi '$ is safe for $\phi$;

3) $\psi = \psi ' \cap \psi ''$, where either $\psi ' $ or $\psi ''$ is safe
for $\phi$.

Let $Sub(\phi)$ denote the set of all subexpressions of $\phi$.
We say that an expression $\phi$ is  safe if 

1)$\phi$ is safe for $\phi$;

2)for every subexpression
$R_\lambda(\psi)$ of $\phi$, $\psi$ is safe for $\phi$.
\end{definition}

%This definition has nothing common with "safe expressions" in computer science.

Some examples of safe expressions  are $x_i$, $R(x)$, $R(R(x)\cap R^{-1}R(x))$,
$$R(R(x)\cap R^{-1}(\top)),\quad R\left(\left(R(x)\cap R^{-1} R(x)\right)\cap\left(R^{-1}(x)\cap R^{-1}(R(x))\right)\right).$$

The Figure \ref{tr2}  shows the dependency tree of the latter expression.

\begin{figure}[tr2]{}
\scalebox{0.9}{\includegraphics{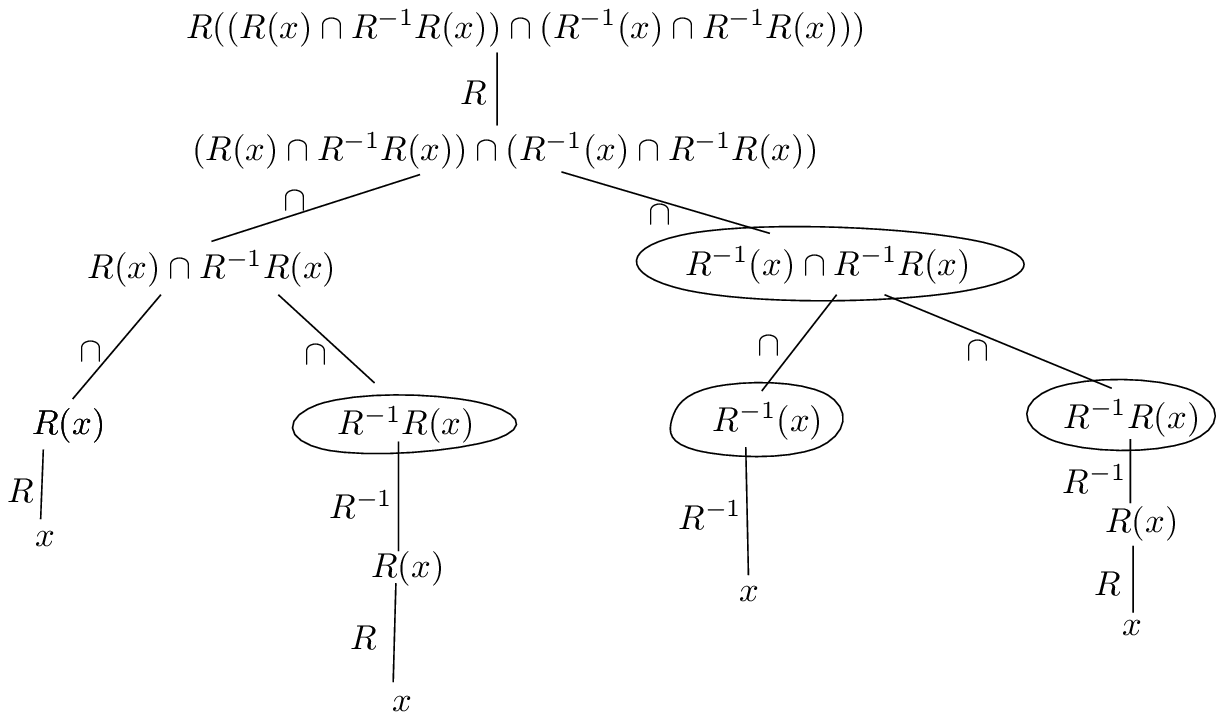}}
\end{figure}

One can  easily check that this expression is safe. 
Denote it by $\phi$. In fact,  
all subexpressions on the left branch are safe for $\phi$, hence $\phi$ 
is safe for itself. However, some of its subexpressions
are not safe for $\phi$; they are circled in the picture. But the operator
$R$ is applied only to the nodes, that are safe for $\phi$.

Examples of non-safe expressions are $R^{-1}(x)$, $R(R^{-1}(x))$, $R(\top)$.

\begin{lemma}\label{safe=k} For any $L$-expression $\phi$
$$\phi\in \K \iffa \phi \mbox{ is safe}.$$
\end{lemma}

\begin{proof*}
By induction on the length of $\phi$. The base is trivial.

Suppose $\phi$ is of the form $\phi_1 \cup \phi_2$. Then $\phi$ is not safe, and
$\phi \notin \K$. 
The same holds if $\phi=R_\lambda^{-1}(\psi)$ or $\phi=R_\lambda^{\b}(\psi)$.

Suppose $\phi=\phi_1 \cap\phi_2$.

If $\phi$ is safe, then either $\phi_1$ or $\phi_2$ is safe. Without any loss of
generality assume that $\phi_2$ is safe. 
Then by the induction hypothesis $\phi_2\in \K$.
Consider $\phi_1$.  Since all subexpressions of the form $R_\lambda(\psi)$ are safe,
$\phi_1$ is of the form 
$POS(\psi_1,\ldots, \psi_k)$, where all  $\psi_i$ are safe. By the  inductive
hypothesis $\psi_i \in \K$, hence $\phi=POS(\psi_1,\ldots, \psi_k)\cap \phi_2\in \K$.

The other way round, if $\phi_1\cap\phi_2 \in \K$, then $\phi_1$ or $\phi_2$
is in $\K$, so either $\phi_1$ or $\phi_2$ is safe, and the other expression
is a positive combination of safe expressions. So $\phi_1\cap\phi_2$ is safe.

Suppose $\phi=R_\lambda(\psi)$. If $\phi$ is safe, then $\psi$ is safe, so $\psi \in \K$,
hence $\phi \in \K$. The other way round, 
if $\phi \in \K$, then $\psi \in \K$, and $\psi$ is safe;
hence $\phi$ is safe. 
\end{proof*}

\begin{lemma} There is a linear algorithm, which 
for any given L-expression $\phi$
decides, whether $\phi$ is safe (or, according to Lemma \ref{safe=k}, whether 
$\phi$ is in $\K$).
\end{lemma}
\begin{proof*}
We run through the syntactic tree of $\phi$ starting from its leaves and assign 
the value 'safe for $\phi$', or 'not safe for $\phi$' to
every node (that is, to a subexpression of $\phi$) according to Definition
\ref{safeexpr}.
If we see that $R_\lambda$ is applied to a node, which is not safe, we stop and
conclude that $\phi$ is not safe. Otherwise, we look whether $\phi$ is safe for
$\phi$, and return the result.

This algorithm takes time proportional to the number of nodes in
the syntactical tree of the expression, hence, it is linear with respect to 
the length of the given expression $\phi$. 
\end{proof*}

\begin{corollary}\label{safe_subst}
Let $\phi$  and $\psi$ be safe expressions. 
After replacing any occurrence of $x_i$ in $\phi$ with $\psi$
we obtain a safe expression $\phi'$.
\end{corollary}

\begin{lemma}[Soundness of $\K$ with respect to $KVF$]\label{safe_sound} 
Let $\phi$ be a regular box-formula of rank $k$ with a head $p^k_i$, 
and let $A$ be the set of all regular
box-formulas of ranks $\le k$, let  $f$ be a map $\{x_1, \ldots, x_n\} \to \P(A)$.
Then $KVF_f^\phi(x_j)$ is in $\K$, and hence, $KF_f^{p^k_i}$ is a union of elements of
$\K$.
\end{lemma}
\begin{proof*} By induction on the length of $\phi$ within the induction on $k$.
The case $k=0$ is trivial. So suppose $0\le l < k$.

Case 1. Let $\phi = p^k_i$. Then $KV^{\phi}=\#$ and $KVF_f^\phi(x_j)=x_j \in \K$.

Case 2. Let $\phi = POS\to \psi$. Then $KV^{\phi}=KV^{\psi}(\# \cap KP^{POS})$ and
$$KVF_f^\phi(x_j)=[KF_f^{p^l_i} / P^l_i ]KV^{\psi}(x_j \cap KP^{POS})=$$
$$=[ KF_f^{p^l_i} / P^l_i ]KV^{\psi}(x_j \cap [
KF_f^{p^l_i} / P^l_i ]KP^{POS})\in \K.$$
In fact, 
$[KF_f^{p^l_i} / P^l_i ]KV^{\psi}(x_j)=KVF^{\psi}_f(x_j)\in \K$ by 
the induction
hypothesis. Consider  $x_j \cap [
KF_f^{p^l_i} / P^l_i ]KP^{POS}$.
Since $l<k$, by the induction hypothesis $KF_f^{p^l_i}$ is a union of safe
expressions. Hence $[KF_f^{p^l_i} / P^l_i ]KP^{POS}$ is a positive
combination of safe expressions.
That is  $x_j \cap [KF_f^{p^l_i} / P^l_i ] KP^{POS} \in \K$
and it is sufficient to apply Corollary \ref{safe_subst}.

Case 3. Let $\phi = \b_\lambda \psi$. Then $KV^{\phi}=KV^{\psi}(R_\lambda(\#))$.
Similarly, $KVF^{\phi}_f(x_j)$ is the result of replacing 
a single occurrence
of $x_j$ with $R_\lambda(x_j)$ in  $KVF^{\psi}_f(x_j)$, which is
safe by the induction hypothesis. So Corollary \ref{safe_subst} implies that $KVF^{\phi}_f(x_j)\in\K.$
\end{proof*}

\begin{lemma}[Completeness of $\K$ with respect to $KVF$]\label{safe_complete}
Let $E(x_1, \ldots, x_k)$ be a safe  $L$-expression and $A$ be the set of
all regular formulas. Then there exists a function
$$f^E: \{x_1, \ldots, x_k\} \to \P(A),$$  and a formula $\phi\in\cup_if^E(x_i)$
with the head $p^l_i$ 
such that $E(x_1, \ldots, x_k)=KF_{f^E}^{p^l_i}=KVF^\phi_{f^E}$. 
\end{lemma}

\begin{proof*}
Induction on the length of $E$. 

 The case $E=x_i$ is
trivial: $$f(x_j)=\left\{\begin{array}{cc}
\emptyset,& \mbox{ if } j\neq i\\
P^0_i,& \mbox{ if } j = i\\
\end{array}\right.$$

Consider an arbitrary safe $E$. Then in the
syntactical tree of $E$ there is a path connecting 
$E$ with some $x_i$, and passing 
only through safe subexpressions of $E$. We denote the subexpressions
on this path by $E_0=\{x_i\}, E_1, \ldots, E_b=E$.

Consider the case $E_1 =R_\lambda(x_i)$. Then consider $E'$  obtained
from $E$ by replacing the subexpression $E_1$ with an expression $E_0$ (that
is we replace $R_\lambda(x_i)$ with $x_i$). Now we apply the induction hypothesis
to $E'$ and obtain a function $f^{E'}$, and a formula $\phi$ with
the head $p^l_i$. 
Then we replace $\phi$   by $\b_\lambda\phi$ in $f^{E'}$, leaving
 $p^l_i$ and others components of 
$f^{E'}$ as they are. This yields us  a function $f^E$, since   
$KV^{\b_\lambda\psi} =  KV^\psi(R_\lambda(\#))$ and the substitution, transforming
$KV$ into $KVF$ is the same for $E$ and $E'$.

Now consider the case $E_1 = \{x_t\} \cap POS(\psi_1,\ldots, \psi_k),$ 
where all $\psi_j$ are safe. By the induction
hypothesis, for any  $\psi_j$ there exist functions $f^{\psi_j}$ and variables $p_j^{l_j}$. 
Let $E'$ be an expression obtained from $E$ by replacing 
$E_1$ with $E_0$. By the induction hypothesis there exists a function
$f^{E'}$ and a formula $\phi$ with the head $p^l_m$. 
Without any loss of generality we may assume that the functions $f^{\psi_j}$ and $f^{E'}$ do not have common propositional
variables and that $l>l_j$ for all $j$ from 1 to $k$. 
Take $f^{E'}$, and replace $\phi$ by $POS' \to \phi$, where $POS'$ is obtained from 
an  L-expression POS
by replacing each of subexpressions $\psi_j$ with $p_j^{l_j}$, $\lor$ with $\cup$,
$\land$ with $\cap$,  $R^{-1}_\lambda$ with $\d_\lambda$, $R^{\b}_\lambda $ by
$\b_\lambda$. We denote the result by
$f^{E'}_{'}$.
Put 
$f^{E}(x_i) =f^{\psi_1}(x_i)\cup \ldots \cup f^{\psi_k}(x_i)\cup f^{E'}_{'}(x_i)$
and the variable $p^l_m$. 
\end{proof*}

\begin{corollary}\label{c18}
Let $\E$ be a set of safe expressions. Then there exists a function 
$f^\E: \{x_0, \ldots, x_k\} \to \P(A)$ and the collection
of variables $\{p^E| E\in\E\}$
\footnote{According to our notation, $p_E$ actually denotes $p^{l_E}_{i_E}$} 
such that for all $E\in\E$ $$E=KF_{f^\E}^{p_E}.$$
\end{corollary}  

\begin{proof*}
By Lemma \ref{safe_complete} for  each $E\in \E$ there exist $f^E$  and $p_E$
such that $E=KF_{f^E}^{p_E}$. Without
any loss of generality we may assume that for different $E$ $f^E$ do not have
common propositional variables. Then we can put
$$f^{\E}(x_i) =\bigcup_{E\in\E} f^E(x_i).$$  
\end{proof*}

Now we see that the class $\K$ describes the values of $KVF$. So 
the values of $KF$ are 
in the closure of $\K$ under $\cup$.

\begin{remark}
This definition of safety does not coincide with the
notion of `safety under bisimulations' from \cite{Bl}.

\begin{definitionCite}{Bl}
A first-order formula $\alpha(x,y)$ is called {\it safe under bisimulation} 
if for all Kripke
models $M$ and  $M'$, bisimulation $Z$ between them 
and points $x_0\in M$, $x_0'\in M$ such that $x Z x'$ 
for all $y_0$  if $M\models \alpha  [x\setminus x_0, y\setminus y_0]$
then there is $y_0'\in M'$ such that 
$M'\models\alpha  [x\setminus x_0', y\setminus y_0']$ and $yZy'$.
\end{definitionCite}

One can  generalize this definition to the following.

\begin{definition}
A first-order formula $\alpha(x_1, \ldots, x_n, y)$ is called {\it safe
under bisimulations}
if for for all Kripke
models $M$ and  $M'$, bisimulation $Z$ between them 
and points $x_i\in M$, $x_i'\in M$  ( $1\le i\le n$) such that $x_i Z x'_i$ 
for all $y_0$  if $M\models \alpha  [x_i\setminus x_i^0, y\setminus y^0]$
then there is $y'^{0} \in M'$ such that 
$M'\models\alpha  [x_i\setminus (x_i^0)', y\setminus y_0']$ and $yZy'$.
\end{definition}

We may conjecture that these two definitions of safety (the syntactic safety
from this paper and safety under bisimulation) coincide.
However, this is not the case. Indeed, the formula $y\in R(R(x_1)\cap R(x_2))$
is safe according to our definition, but  not safe under bisimulations.
\end{remark}

\section{Generalized Sahlqvist formulas}

\begin{definitionCite}{Gor2} 
A { \sl generalized Sahlqvist implication } 
is a formula $GSA \to \bot$, where
$GSA$\footnote{Generalized Sahlqvist Antecedent} is  built from regular box-formulas and negative formulas (that is,
negations of positive formulas) 
using only $\land, \lor, \d_\lambda$.
If we prohibit the use $\lor$ in $GSA$, we obtain the definition of a 
{\sl generalized  simple Sahlqvist implication}.

A {\sl generalized Sahlqvist formula\footnote{In subsequent 
publications Goranko and Vakarelov refer to these formulas as
the monadic inductive formulas}} is a formula built up from
generalized Sahlqvist implications by applying boxes and conjunctions, and 
by applying disjunctions only to formulas without common proposition letters.
\end{definitionCite}

The reduction of a generalized Sahlqvist formula to a generalized simple Sahlqvist implication is standard \cite{Bl}.
So without any loss of generality we may 
consider a generalized simple Sahlqvist implication $GSA \to \bot$,
where $GSA$ is built from regular box-formulas and negative formulas using only $\land$ and $\d_\lambda$. It is convenient to represent such formulas with labelled trees of a special kind, similar to syntactical trees.

\begin{definition} Consider a structure
$\hat T=(W, (R_\lambda:\lambda\in\Lambda))$.
{\sl A path from $x_1$ to $x_n$ } in $\hat T$ is a sequence $x_1 \lambda_1 x_2 \lambda_2 x_3 
\ldots x_n$, where $x_i \in W$, $\lambda_i \in \Lambda$ and 
$x_i R_{\lambda_i} x_{i+1}$ in $\hat T$. Two paths $x_1 \lambda_1 x_2 \lambda_2 x_3 
\ldots x_n$ and $x'_1 \lambda'_1 x'_2 \lambda'_2 x'_3 
\ldots x'_n$ are called {\sl equal} if for all $1\le i \le n$
$x_i=x'_i$ and for all  $1\le i \le n-1$ $\lambda_i = \lambda'_i$.

A pair $(\hat T,r)$ is called a {\sl tree with a root $r$} if the
following holds
 
1) $r \in W$,

2) $R_\lambda^{-1} (r) = \emptyset$ for all $\lambda \in \Lambda$,

3) for all  $x\neq r$ there is a unique path from $r$ to $x$. 

Let $A$ be a  set of modal formulas.
A {\sl labelled tree with a root $r$} is a tuple $T=(W, (R_\lambda:\lambda\in\Lambda), r, f)$, where  
$(W, (R_\lambda:\lambda\in\Lambda), r)$ is a tree with a root $r$ and
$f$ (a label function)
is a map from  $W$ to $\P(A)$. 
\end{definition}

\begin{definition}
Let $\phi$ be built up from formulas of $A$ by applying only diamonds
and conjunction. 
A {\sl reduced syntactical tree} of a formula  $\phi$ 
is a labelled tree defined by induction on the length of $\phi$.

Case 1: $\phi = a$, where $a \in A$.
Then $T^\phi$ contains a single  
point $x$. The map $f^\phi$ takes $x$ to $\{ a\}$ and the
relations $R_\lambda^\phi$ are empty.

Case 2: $\phi = \chi \land \psi$. Then put $W^\phi = 
(W^\chi \backslash \{r^\chi\}) \cup (W^\psi\backslash \{r^\psi\}) \cup\{ r^\phi\}$,
where $r^\phi$ is some new point.
The relations $R_\lambda$ on $W^\chi$ and $W^\psi$ remain the same, and
$r^\phi R_\lambda w$ iff $w\in W_\chi$ and $r^\chi R^\chi_\lambda w$
or $w\in W_\psi$ and $r^\psi R^\psi_\lambda w$. The map $f^\phi$ sends $r^\phi$ to 
$f^\chi(r^\chi) \cup f^\psi(r^\psi)$ and is equal to $f^\chi$ or  $f^\psi$ in
all other points.

Case 3: $\phi = \d_\lambda \psi$. Then $W^\phi$ = $W^\psi\cup \{ r^\phi\}$, where $r^\phi$ 
is a new point. The $R_\mu$ for $\mu \neq \lambda$ we leave untouched, and 
to $R_\lambda$ we add an arrow, joining $r^\phi$ with $r^\psi$. We put $f(r^\phi) = \emptyset$,
and do not change $f$ in all  other points.
\end{definition}

\begin{example}
The reduced syntactical tree of the formula 
$$\phi=\d(\b p \land\b q\land \d (\d \b q \land \d \b\b p))\land p$$
is shown in the Figure \ref{tr4}.
\begin{figure}[tr4]{}
\includegraphics{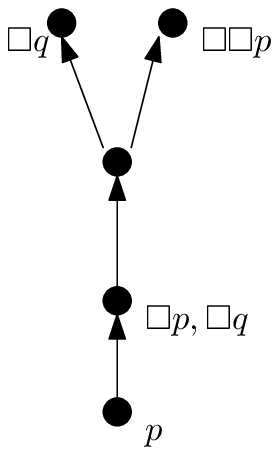}
\end{figure}
\end{example}

\begin{lemma}\label{tree}
 Let $A$ be an arbitrary set of modal formulas and let $\phi$ be built 
from formulas of $A$ using only $\land$ and $\d_\lambda$.
Let $T^\phi=(W^\phi, (R^\phi_\lambda:\lambda \in \Lambda), r^\phi, f^\phi)$ be a reduced syntactical tree of $\phi$. 
Then for all frames $F=(W,R_\lambda :\lambda\in\Lambda)$ for
any valuation $\theta$ $F, x, \theta \models \phi$ iff there exists a 
monotonic map $h: T^\phi \to F$ 
(that is for all $x,y \in W^\phi$ if $x R^\phi_\lambda y$ then 
$h(x) R_\lambda h(y)$) such 
that $h(r^\phi)=x$ and for any $w\in W^\phi$,  $a\in A$ if $a\in f^\phi(w)$ then $F,h(w),\theta \models a$.
\end{lemma}

The proof of the Lemma \ref{tree} trivially follows from the semantics  of $\land$ and
$\d_\lambda$.

For  Sahlqvist formulas $A$ is the set of all boxed atoms and negative formulas. 
For generalized Sahlqvist formulas 
$A$ is the set of all regular box-formulas and negative formulas. 

The next lemma shows the standard second-order quantifier elimination in 
a simple generalized Sahlqvist implication.

\begin{lemma} (cf. \cite{Gor2} and \cite{Bl}, Section 3.6)
Let $\phi$ be a simple generalized Sahlqvist implication $\phi$ with 
a reduced syntactical tree $T=(\{y_0, y_1, \ldots, y_n\}, 
(R^T_\lambda : \lambda \in \Lambda), y_0, f)$. Let 
$f_{REG}(y_i)=f(y_i)\cap REG$, and $f_{NEG}(y_i)=f(y_i)\cap NEG$ where $REG$ and $NEG$ are respectively the sets of  all regular box-formulas and all negative formulas.
Then the first-order correspondent of $\phi$ is of the form
\begin{equation}\label{ssa}
[(x_j\in KF_{f_{REG}}^{p^k_i})^{\#}/P^k_i(x_j)]\forall x_1 \ldots \forall x_n  
\left(\bigwedge_{y_i R^T_\lambda y_j}x_i R_\lambda x_j \to \right.\end{equation}
$$
\left. \bigvee_{\psi\in f_{NEG}(y_j)} (x_j \models\neg \psi)^*
\phantom{\bigvee_{y_i R^T y_j}x_i}\right).$$
Here $KF_{f_{REG}}^{p^k_i}$ is the minimal valuation 
( see Definition \ref{K_f}),
$^{\#}$ denotes the first-order transcription of 
$x_j \in KF_{f_{REG}}^{p^k_i}$, defined
on the page \pageref{diez}, and $^*$ means the standard first-order translation
of a modal formula.
\end{lemma} 

\begin{proof*}
The proof is standard. As in  the Sahlqvist theorem, we can eliminate
the second-order quantifiers by substituting appropriate (minimal) valuations.

Let $\phi$ be a simple generalized Sahlqvist implication with
a reduced syntactical tree 
$T=(\{y_0,y_1, \ldots, y_n\}, (R^T_\lambda:\lambda \in \Lambda),y_0, f)$. 
Then for any frame  $F$,  
$F, x_0 \models\phi$ is equivalent to  the universal
second order formula
$$\forall P^{k_1}_{i_1} \ldots \forall P^{k_m}_{i_m} 
\left(\exists x_1 \ldots \exists x_n 
\left(\bigwedge_{y_i R^T_\lambda y_j}x_i R_\lambda x_j \land\bigwedge_{i} 
(x_i \models f(y_i))^*\right) \to \bot \right),$$
where for a set of modal formulas $f(y_i)$ the notation $x_i \models f(y_i)$ 
means that in the point $x_i$ all members of $f(y_i)$ are true.

Now we can put the existential quantifiers in the prefix. Since
they are in the antecedent of the implication, they become
universal:
$$\forall P^{k_1}_{i_1} \ldots \forall P^{k_m}_{i_m} \forall x_1 \ldots \forall x_n\left(  
\left(\bigwedge_{y_i R^T_\lambda y_j}x_i R_\lambda x_j \land\bigwedge_{i} (x_i \models f(y_i))^*\right) \to \bot\right)$$ 
Then let us swap them with the second-order quantifiers:
$$\forall x_1 \ldots \forall x_n\forall P^{k_1}_{i_1} \ldots \forall P^{k_m}_{i_m} \left(  
\left(\bigwedge_{y_i R^T_\lambda y_j}x_i R_\lambda x_j \land\bigwedge_{i} (x_i \models f(y_i))^*\right) \to \bot \right).$$

Now we apply the equivalence 
$A\land B \to C\equiv A\to(B\to C)$, yielding
$$\forall x_1 \ldots \forall x_n\forall P^{k_1}_{i_1} \ldots \forall P^{k_m}_{i_m} \left(  
\bigwedge_{y_i R^T_\lambda y_j}x_i R_\lambda x_j \to
\left(\bigwedge_{i} (x_i \models f(y_i))^* \to \bot\right) \right).$$
Let us move the second-order universal quantifiers to the consequent:
$$\forall x_1 \ldots \forall x_n  
\left(\bigwedge_{y_i R^T y_j}x_i R x_j \to  \forall P^{k_1}_{i_1} \ldots \forall P^{k_m}_{i_m} \left( \left(
\bigwedge_{i} (x_i \models f(y_i))^*\right) \to \bot\right)\right)$$ 

Now let us recall that $f(y_i)= f_{NEG}(y_i) \cup f_{REG}(y_i)$.
Let us  move the formulas of $f_{NEG}(y_i)$ from the antecedent to the consequent of
the inner implication:  
$$\forall x_1 \ldots \forall x_n  
\left(\bigwedge_{y_i R^T y_j}x_i R x_j \to  \forall P^{k_1}_{i_1} \ldots \forall P^{k_m}_{i_m} \left( \left(
\bigwedge_{i}\bigwedge_{\psi\in f_{REG}(y_i)} (x_i \models \psi)^*\right) \to
\right.\right.
$$$$ \left.\left.\bigvee_{i}\bigvee_{\psi\in f_{NEG}(y_i)} \neg(x_i \models \psi)^*)\right)\right).$$
According to Lemma \ref{theta_min}, there is the smallest valuation
verifying the antecedent 
$$\bigwedge_{i}\bigwedge_{\psi\in f_{REG}(y_i)} (x_i \models \psi)^*.$$
 The negation of a negative formula is positive. So we can eliminate the second-order quantifiers by substituting the minimal valuation  and obtain the formula (\ref{ssa}).
\end{proof*}

\section{Generalized Kracht Formulas}

Now we will extend Kracht's theorem to generalized Sahlqvist 
formulas. To this end we need an extension of our first-order language.
The only contribution of this work is the usage of $\K$. 
All other definitions from this section 
(restricted quantification, inherently universality) are taken from \cite{Bl} 
and originate from Kracht.

We abbreviate the first-order formula $\forall y (xR_\lambda y \to \alpha(y))$ to
$(\forall y \tr_\lambda x) \alpha (y)$.
Likewise $\exists y (x R_jy \land \alpha(y))$ is
abbreviated to $(\exists y\tr_\lambda x)\alpha(x)$. 
%We call the constructs 
%$(\forall y \tr_\lambda x)$ and $(\exists y \tr_\lambda x)$ {\sl restricted %quantifiers}.
We shall use only formulas, in which variables
do not occur both as free and bound, and in which two distinct occurrences of quantifiers do not bind
the same variable; we  call such formulas {\sl clean}.

%Let $\bar\K=\K\cup\{\top, \bot\}$, where $\K$ is the class of all
%safe expressions from Section 3.
Let $\K$ be the class of all
safe expressions from Section 3.
We  add new $(k+1)$-ary predicates  
$x_l\in E(x_1, \ldots, x_k)$ for
any  expression $E\in \K$. Depending on the context, they can  also be considered
as abbreviations for the corresponding first-order formulas with free variables
$x, x_1, \ldots, x_k$.

More precisely, for any  $L$-expression $E$ (not necessary safe) 
we define a first-order formula
$(x_l \in E)^\#$\label{diez} by the recursion on the length of $E$:
$$(x_l \in x_i)^\#:= (x_l = x_i);$$ 
$$(x_l \in \top)^\#:= (x_l = x_l);$$ 
$$(x_l \in \bot)^\#:= \neg(x_l = x_l);$$ 
$$(x_l \in E_1\cap E_2)^\#:= (x_l \in E_1)^\# \land (x_l \in E_2)^\#;$$ 
$$(x_l \in E_1\cup E_2)^\#:= (x_l \in E_1)^\# \lor (x_l \in E_2)^\#;$$ 
$$(x_l \in R_\lambda^{-1}(E))^\#:= 
\exists y (x_l R_\lambda y \land (y \in E)^\#);$$ 
$$(x_l \in R_\lambda^{\b}(E))^\#:= 
\forall y (x_l R_\lambda y \to (y \in E)^\#);$$ 
$$(x_l \in R_\lambda(E))^\#:= 
\exists y (y R_\lambda x_l \land (y \in E)^\#).$$
This translation obviously corresponds to the standard set-theoretic
semantics.

\begin{definition} (cf \cite{Bl}, p. 172)
We call a formula {\sl restrictedly positive} if it is built up from  
formulas $y\in E(x_1, \ldots, x_k)$, using 
$\land$, $\lor$ and restricted quantifiers.

We say that a variable $x$ in a clean formula $\alpha$ is
{\sl inherently universal} if either $x$ is free, or  $x$ is bound by a restricted 
universal quantifier which is not
in the scope of an existential quantifier. 

A formula $\alpha$ in the 
extended first-order language is
called a {\sl generalized Kracht formula with free variables} if 
$\alpha$ is clean, restrictedly positive and  
in every subformula of the form $y\in E(v_1, \ldots, v_k)$ ($E\in \K$),  
the variables 
$v_1, \ldots, v_k$ are inherently 
universal.  A formula $\alpha$ is called 
a {\sl generalized Kracht formula } if it is a generalized Kracht formula with free variables and it contains exactly one free variable.
\end{definition}

The definition of ordinary Kracht formulas is obtained from this definition
by replacing $\K$ with $\{R_{\lambda_1}\dots R_{\lambda_n}(x_j)\}\cup
\{R_{\lambda_1}^{-1}\dots R^{-1}_{\lambda_n}(x_j)\}$.

Now we are ready to state the main theorem.

\begin{theorem} \label{main}
A first-order formula $\phi$  is a first-order
correspondent of a generalized Sahlqvist formula iff $\phi$ is a generalized Kracht 
formula.
\end{theorem}

%Before proving the theorem, we do some comments on this definition and give some %examples of such formulas.

Note, that every ordinary Kracht formula can be rewritten as a 
generalized Kracht formula.
Namely, instead of $x R_{\lambda_1} \ldots R_{\lambda_k} y$, 
where $x$ is inherently universal, we
write $y\in R_{\lambda_k} \ldots R_{\lambda_1}(x)$ 
(obviously, $R_{\lambda_k} \ldots R_{\lambda_1}(x)$ is a safe expression).
Instead of $y R_{\lambda_k} \ldots R_{\lambda_1} x$, 
where $x$ is inherently universal, we write
$$(\exists z_1\tr_{\lambda_1} y)(\exists z_2\tr_{\lambda_2}z_1) \ldots
(\exists z_k\tr_{\lambda_k} z_{k-1}) (z_k \in  x).$$
%The converse is not true, as next example shows; the
%details can be found in \cite{Gor2}. 

\begin{example}\label{example}
Consider the formula $$D_2=p\land \b(\d p \to \b q) \to \d\b\b q$$ 
from \cite{Gor2}.
Its first-order correspondent is a generalized Kracht formula
$$FO(D_2)=\exists y \tr x\left( \forall z'\tr y \forall z\tr z' 
z \in R(R(x)\cap R^{-1}(x))\right),$$ 
or, in a more standard form, 
$$FO(D_2)=\exists y \left(x R y \land \forall z \left(y R^2 z \to 
z \in R(R(x)\cap R^{-1}(x))\right)\right).$$

In \cite{Gor2} the authors show that it is not equivalent to any standard 
Sahlqvist formula.

Consider the formula
$$ ns=p\land \b_1(\d_1 p \to  \b_3 r) \to  \d_2 (\d_2 p\land \d_3 r).$$
Then
%$$FO(ns)=\exists y \left(x R_1 y \land y\in R_1^{-1}(x)\land\exists v \left(y R_3 v %\land 
%v \in R_3(R_2(x)\cap R_2^{-1}(x))\right)\right).$$
$$FO(ns)=\exists y \tr_1 x \left(y\in R_1^{-1}(x)\land\exists v \tr_3 y \left( 
v \in R_3(R_2(x)\cap R_2^{-1}(x))\right)\right).$$
This generalized Kracht formula 
is equivalent to
$$\exists y \exists z \exists v
(xR_1y\land yR_1x\land xR_2z
\land zR_2x\land yR_3v \land z R_3 v).$$

The formula $cub_1$ is a theorem  of  $\mbox{\bf K}^3$ \cite{Sh78}, see also
\cite{GZ}, p.~397
$$\label{cub1}
\begin{array}c
cub_1 = \left[ \d_1(\b_2p_{12}\land \b_3p_{13})\land \d_2(\b_1p_{21}\land
\b_3p_{23})\land \d_3(\b_1p_{31}\land \b_2p_{32})\land\right.\\
\b_1\b_2(p_{12}\land p_{21}\kern -0.25em \to\kern -0.2em \b_3 q_3)
\kern -1pt \land\kern-1pt
\b_1\b_3(p_{13}\land p_{31}\kern -0.25em\to\kern -0.2em \b_2 q_2)\kern-1pt\land\kern-1pt
\b_2\b_3(p_{23}\land p_{32}\kern -0.25em\to\kern -0.2em \b_1 q_1)\left.\right]\\
\to \d_1\d_2\d_3 (q_1\land q_2\land q_3).
\end{array}
$$

Its first-order correspondent is a generalized Kracht formula
$$\forall x_1 \tr_1 x\forall x_2 \tr_2 x\forall x_3\tr_3 x \exists y'\tr_1 x
\exists y''\tr_2 y' \exists y\tr_3 y''$$$$(y \in 
R_3(R_2(x_1)\cap R_1(x_2))\land y\in  R_2(R_3(x_1)\cap R_1(x_3))\land$$$$
\land y \in R_1(R_2(x_3)\cap R_3(x_2))). $$
%$$\forall x_1 \tr_1 x\forall x_2 \tr_2 x\forall x_3\tr_3 x \exists y
%((x R_1 R_2 R_3 y)\land$$$$\land y \in 
%R_3(R_2(x_1)\cap R_1(x_2))\land y\in  R_2(R_3(x_1)\cap R_1(x_3))\land$$$$
%\land y \in R_1(R_2(x_3)\cap R_3(x_2))). $$
This formula is equivalent to (\ref{cub_fo}).
 \end{example}

Examples of generalized Kracht formulas applied to many-dimensional
modal logics can be found in \cite{Agi1},
\cite{Agi2} and \cite{K3}.

The rest of the paper will be devoted to the proof of this theorem.

\section{Quasi-safe expressions}
\begin{definition}
An $L$-expression is called {\sl quasi-safe } if it is a positive combination
of safe expressions.
\end{definition}

The expression $\top$, $\bot$, $R^{-1}(\top)$ are here considered as 
quasi-safe but not safe.

%It follows directly from definition of $\K$ that an intersection of a quasi-safe
%expression and a safe expression is a safe expression.

If we extend our first-order 
language with atomic formulas $x\in E$ where $E$ is a quasi-safe
expression, we obtain a quantifier elimination in the scope of 
the existential quantifier.

\begin{lemma}\label{qe}
Let $\psi$ be a generalized Kracht formula with free variables, such that all
atomic formulas of $\phi$ are of  the form $y\in E$ where all variables occuring 
in $E$ are free. 
Then $\psi$ is equivalent to a quantifier free formula $\psi'$ in the
extended language 
(cf. \cite{Bl}, p. 175). 
\end{lemma}
\begin{proof*}
We apply the 
induction on the number of quantifiers in $\psi$.

Consider the case $\psi=\exists y \tr_\lambda x \phi$. By the induction hypothesis,
$\phi$ is a quantifier free formula. Hence we can assume that it is of the form
$\phi = K_1 \lor \ldots \lor K_n$, where $K_i$ are conjunctions of atomic formulas.
But then $\psi\equiv \exists y\tr_\lambda x K_1 \lor\ldots \lor \exists
y\tr_\lambda  x K_n$. Then, since all $E_i$ do not contain $y$, 
we can transform each of the disjuncts as follows
$$\exists y\tr_\lambda x (\alpha_1 \in E_1 \land \ldots \land \alpha_m \in E_m)
\equiv  \bigwedge_{\alpha_i \neq y } \alpha_i \in E_i \land
x\in R_\lambda^{-1}\left( \bigcap_{\alpha_i = y} E_i\right),$$
and obtain a quantifier free equivalent of $\psi$.

Similarly, let  $\psi=\forall y \tr_\lambda x\phi$. By the induction hypothesis, 
$\phi$ is  quantifier free,  so it can be presented in the form
$\phi = D_1 \land \ldots \land D_n$, where $D_i$ are disjunctions 
of atomic formulas.
But then $\psi$ is equivalent to 
$\forall y\tr_\lambda x D_1 \land\ldots \land \forall
y\tr_\lambda  x D_n$. Then each of conjucts can be transformed as follows
$$\forall y\tr_\lambda x (\alpha_1 \in E_1 \lor \ldots \lor \alpha_m \in E_m)\equiv 
\bigvee_{\alpha_i \neq y} \alpha_i \in E_i \lor
x\in R_\lambda^{\b}\left( \bigcup_{\alpha_i = y} E_i\right).$$
\end{proof*}

\begin{corollary}\label{qe1}
Let $\psi$ be a generalized Kracht formula, beginning with an existential quantifier.
Then $\psi$ is equivalent to a quantifier free formula $\psi'$ in the
language with quasi-safe atoms.
\end{corollary}

\section{Proof of the theorem.}
`Only if'. If $\phi$ is a simple generalized Sahlqvist
implication, then the statement follows from (\ref{ssa}). 
It is sufficient to
note that 
$$\forall x_1 \ldots \forall x_n  
\left(\bigwedge_{y_i R^T_\lambda y_j}x_i R_\lambda x_j \to C\right)$$ is
equivalent to 
$$\forall x_1 \tr_\lambda x_{p(1)}\ldots\forall x_n\tr_\lambda x_{p(n)} C,$$
where $y_{p(i)}$ is the unique predecessor  of $y_i$ in $T$. The variables
$x_1, \ldots, x_n$ are inherently universal, the 
disjunction $C$ is built from atomic formulas using $\lor, \land$ and restricted
quantifiers,  and every atomic formula is
of the form $v\in E(x_{i_1}, \ldots, x_{i_k})$, since 
we substitute the disjunctions of such formulas 
for all $P^k_i$ in the standard translation of positive formulas.

The general case follows from Lemma 3.53 of \cite{Bl} stating  that
\begin{itemize}
\item if $\phi$  and $\alpha(x)$ are locally correspondents, so are 
$\b_\lambda\phi$ and
$\forall y\tr_\lambda x \alpha(y)$,
\item if $\phi$ locally corresponds to $\alpha(x)$ and $\psi$ locally corresponds to $\beta(x)$ then $\phi\land \psi$ locally corresponds to $\alpha(x) \land\beta(x)$,
\item if $\phi$ locally corresponds to $\alpha$, $\psi$ locally corresponds to 
$\beta(x)$ and $\phi$ and $\psi$ do not have propositional letters in common, then 
$\phi \lor \psi$ locally corresponds to $\alpha(x) \lor \beta(x)$,
\end{itemize}
and it remains to note that the class of generalized Kracht formulas is closed under disjunction, conjuntion and necessitation.

To prove `if', 
we need to  generalize the notion of modal definability to
first-order formulas with  many free variables 
(cf. \cite{Kr2}, p. 193).

We say that a first-order formula $\phi(x_1, \ldots, x_n)$ is 
{\sl definable } if 
there is a sequence of modal formulas $\phi_1, \ldots, \phi_n$ such that
for any frame $F=(W, (R_\lambda:\lambda\in\Lambda))$  for any points $x_1^0, \ldots, x_n^0 \in W$
$$F\models \Phi(x_1, \ldots, x_n) [x_1^0, \ldots x_n^0] \iffa
\begin{array}l
\mbox{ for any valuation $\theta$ there exists $i$}\\
\mbox{ such that } F,x_i^0, \theta \models \phi_i.
\end{array}$$

Here the left hand $\models$ means the truth in $F$ considered as a 
classical first-order structure.

For example, a formula $x_1 R x_2 $ is definable by the sequence 
$\d\neg p, p$. Clearly, that if $\phi$ has a single variable, 
then this definition coinsides with the standard modal definability.
 
Now we show that  the formula
$(x_l \in E)^\#$ is definable for all quasi-safe $E$.

%The formula  $x_l \in E(x_1, \ldots, x_k)$ where $E$ is a safe expression
%is describable too. Indeed, by \ref{safe_complete} there exists 
%$f: \{x_1, \ldots, x_k\} \to \Ml$ and a variable $p$, such that 
%$E=KF^p_f$. So, $x_l \in E(x_1, \ldots, x_k)$ is describable by a cortege
%$\neg f(x_1), \ldots, \neg f(x_k), p$. Moreover, the formula
%$x_l \in E(x_1, \ldots, x_k)$ is describable even if $E$ is quasi-safe.

To this end, consider the following translation $^T$ 
from quasi-safe expressions to modal language.
Let $E$ be a quasi-safe expression. Let $\E$ be the set of all safe subexpressions
occuring in $E$.

Now we define $E^T$ by the induction on the length of   $E$:

if  $E$ is safe then $E^T= p_E$;

if $E=E_1 \cap E_2$ then $E^T=E_1^T \land E_2^T$;

if $E=E_1 \cup E_2$ then $E^T=E_1^T \lor E_2^T$;

if $E=R^{-1}_\lambda E_1$ then $E^T = \d_\lambda E_1^T$;

if $E=R^{\b}_\lambda E_1$ then $E^T = \b_\lambda E_1^T$.

\begin{lemma}\label{fE}
Let $E$ be quasi-safe and let  $\E$ be the set of all safe subexpressions
occuring in $E$. 
Let $f^\E$ be 
the function from Corollary \ref{c18} for the set $\E$. 
Then $(x_l\in E)^\#$ is definable
by the sequence $\bar \phi = \phi_1, \phi_2, \ldots, \phi_m$, such that
$$\phi_i = \left\{
\begin{array}{ll}
\bigvee_{\phi\in f^\E(x_i)} \neg \phi, & i\neq l;\\
\bigvee_{\phi\in f^\E(x_i)} \neg \phi\lor E^T, & i= l.\\
\end{array}\right.
$$
\end{lemma}
\begin{proof*}
Suppose that we have a frame $F=(W, (R_\lambda:\lambda))$, and the variables
$x_i$ are identified with points of $W$. Then we can evaluate $E$ and regard
it as a subset of $W$.

Let us call a valuation $\theta$ {\sl admissible }  
if for all $i$ $x_i,\theta \models f^{\E}(x_i)$.

Consider the following statements:

$(1)$ $x_l \in E$

$(2)$ $x_l, \theta_{\min} \models E^T$, where $\theta_{\min}$ is the valuation from
Lemma \ref{theta_min}.

$(3)$ for all admissible valuations $\theta,  x_l \models E^T$.

Then due to the form of $\phi_i$, 
the statement of the lemma can be rephrased as $(1) \iffa (3)$. But
Lemma \ref{theta_min} ensures that  $(2)\iffa (3)$. 

Let us prove $(1) \iffa (2)$ 
by induction on the length of a quasi-safe $E$.
 
The base. Suppose $E$ is safe. In this case  $E^T=p_E$ and  
$$\theta_{\min}(p_{E})=KF_{f^\E}^{p_E}=E.$$
The first equality holds by Lemma \ref{theta_min} and
the second one by Corollary \ref{c18}, and 
the statement  is clear.

The induction step trivially follows from the interpretation of
$\lor, \land, \d_i,$ and  $\b_i$ in Kripke semantics.

In fact, 
let $E=E_1 \cap E_2,$ that is  $E^T=E_1^T\land E^T_2$.
$$x_l\in E_1\cap E_2 \iffa 
x_l \in E_1 \mbox{ and } x_l \in E_2 \iffa$$$$\iffa
 \theta_{\min},  x_l \models E_1^T\mbox{ and } 
 \theta_{\min}, x_l \models E^T_2\iffa
\theta_{\min}, x_l \models E_1^T\land E^T_2.
$$
The case of the disjunction is similar.

Let $E=R^{-1}_\lambda(E_1)$. Then 
$$x_l \in R^{-1}_\lambda(E_1) \iffa \exists y
(x_l R_\lambda y \land (y\in E_1)) \iffa$$$$
\iffa
\exists y (x_l R_\lambda y \mbox{ and } \theta_{\min}, y \models E_1^T ) \iffa
\theta_{\min}, x_l \models \d_\lambda E_1^T.$$

Let $E=R^{\b}_\lambda(E_1)$. Then 
$$x_l \in R^{\b}_\lambda(E_1) \iffa \forall y
(x_l R_\lambda y \to  (y\in E_1)) \iffa$$$$
\iffa
\forall y (x_l R_\lambda y \to  \theta_{\min}, y \models E_1^T)  \iffa
\theta_{\min}, x_l \models \b_\lambda E_1^T.$$
\end{proof*}

We also need a dual version
of  Theorem 5.6.4 from \cite{Kr2}:
\begin{theoremCite}{Kr2} \label{Kr}
If $\alpha(x_0)$ is obtained from definable formulas using
conjunction, disjunction  and restricted universal quantification, then
$\alpha(x_0)$ is definable.
\end{theoremCite}

Now we are ready to prove the main theorem.

\begin{lemma}\label{last}
Let $\alpha(x_0)$ be a first-order formula with the only free variable $x_0$. 
Then the  following statements are equivalent:

(1) $\alpha(x_0)$ is a first-order correspondent of a generalized Sahlqvist formula;

(2) $\alpha(x_0)$ is a generalized Kracht formula;

(3) $\alpha(x_0)$ is obtained from formulas of the form $x_l\in E$, 
where $E$ is quasi-safe, using conjuction, disjunction and  restricted universal quantification.
\end{lemma}
\begin{proof*}

$(1) \to (2)$  was proved at the beginning of Section 7, in the 
'only if' part.

$(2) \to (3)$. Given a generalized Kracht formula $\phi$, we 
apply the quantifier elimination from Corollary \ref{qe1} to its
maximal subformulas beginning with existential quantifiers.  
Then we obtain a formula satisfying (3).

$(3) \to (1)$. Apply Lemma \ref{fE} and Theorem \ref{Kr} to $\alpha(x_0)$.
\end{proof*}

It is clear that Lemma \ref{last} implies Theorem \ref{main}

\section{Discussion}
1. The papers \cite{Gor1}, \cite{Gor2} deal mainly with
`inductive' formulas, that are, in brief, 
generalized Sahlqvist formulas in polyadic modal languages.
The theory of inductive formulas is in some sense more elegant,
than the theory of generalized Sahlqvist formulas. So it would be
interesting to extend Kracht's theorem to inductive formulas
in polyadic modal languages.  D. Vakarelov made a conjecture that 
their characterization may be nicer.

2. Note that there is a certain asymmetry between $R$ and $R^{-1}$ in
the definition of safe expressions. In  temporal language 
this asymmetry disappears, and, as Gorando and Vakarelov
show in \cite{Gor2}, every generalized
Sahlqvist formula is semantically 
equivalent to the standard Sahlqvist one. 

%This result can be obtained as an easy consequence of generalized Kracht theorem.
%In fact, if we have symbols $R$ and $R^*$ bound by the correlation
%$R^* = R^{-1}$, and a

3. Traditionally the correspondence between Sahlqvist and Kracht formulas
and their generalization 
is considered from the viewpoint of definability.
We have several answers  to the natural question 
``what first-order formulas are modally definable?''
 For example there is a  sufficient 
syntactic condition
given by the class of   Kracht formulas and their generalization, 
and there is  
also a semantical characterization given by Goldblatt-Thomason theorem \cite{GT}. 
But we can also ask
when the modal logic of an elementary class is finitely axiomatizable.
Kracht formulas and their generalization 
give a sufficient syntactic condition in this case too,
but we do not have a semantical characterization. 
It would be interesting to look for other elementary classes
with finitely axiomatizable modal logics. For example, it is known 
\cite{BSS} that the modal logic of the elementary class of the formula 
$\exists y (x R y \land  R(y)\subset \{y\})$ is finitely axiomatizable.

\bibliography{template}
\end{document}